# Les autres de l'un :
# deux enquêtes prosopographiques sur Charles Hermite

Catherine GOLDSTEIN

La prosopographie s'attache en principe aux collectivités. Les auteurs de la synthèse *A Short Manual to the Art of Prosopography* opposent ainsi naturellement prosopographie et biographie : « [Prosopography] is rather a research approach than a method *sui generis* ; an attempt to bring together all relevant biographical data of groups of persons in a systematic and stereotypical way. […] Biography studies the particular life histories of individuals. […] Biographies are written mainly about exceptional or special people and try to achieve a better understanding of their personality. Prosopography is not interested in the unique, but in the average, the general and the 'commonness' in the life histories of more or less large number of individuals. The individual and the exceptional is important only insofar as it provides information on the collective and the 'normal'. For a prosopographer, extraordinary people […] are less appealing and to some extent even disturbing because they are out of the ordinary[1]. »

Relever systématiquement, pour un corpus prédéterminé de personnes, les valeurs de certaines caractéristiques (elles-mêmes prédéterminées) doit permettre d'y repérer des constances, des régularités, des effets de masse, des sous-populations intéressantes[2]. Une société, un certain collectif, sont ainsi appréhendés par une juxtaposition et une analyse de micro-informations sur chaque constituant, chaque membre : quels sont par exemple les revenus, les fonctions sociales, le nombre d'enfants, les distractions, le nombre d'heures de travail des ouvriers d'une même fabrique ou des instituteurs à une certaine époque ? Cette problématique contraint les cadres tant pratiques que théoriques de l'enquête, et jusqu'aux critiques mêmes qui ont pu être faites sur leur articulation : la délimitation préalable des catégories du questionnaire ne reproduit-elle pas des préjugés sur la population étudiée ? Peut-on interpréter des informations à l'intérieur d'un groupe, sans tenir compte de ses rapports à d'autres groupes, ou même comparer les entrées d'une même catégorie pour plusieurs membres lorsqu'elles ont été parfois obtenues par l'intermédiaire de sources variées, dépendant des membres considérés ?

Les biographies d'hommes et de femmes célèbres posent a priori d'autres problèmes : la partialité potentielle de la perspective, que des sources centrées autour d'une seule personne tendent à créer, la complexité des liens entre la vie et l'œuvre, et leur existence même, la continuité peut-être artificielle entre les événements induite par la seule chronologie, sont autant de questions classiques du genre. Les prosopographies d'artistes ou de scientifiques, par nature, visent donc à quelque hétérodoxie : il semble que l'on ne puisse traiter des écrivains ou des peintres collectivement si l'on se place dans une conception de l'activité artistique qui est celle de l'œuvre de génie, de l'exceptionnel. L'approche prosopographique a donc particulièrement fleuri sur fond de critique sociale de l'art ou de la science. Le sociologue est-il alors « voué au relativisme, au nivellement des valeurs, à l'abaissement des grandeurs, à l'abolition des différences qui font la singularité du 'créateur', toujours situé du côté de l'Unique ? Cela parce qu'il aurait partie liée avec les grands nombres, la moyenne, le moyen, et, par conséquent, avec le médiocre, le mineur, les *minores,* la masse des petits auteurs obscurs »*,* s'interroge ainsi Pierre Bourdieu à propos des écrivains[3]. La même question se pose pour les scientifiques : David Sturdy rappelle que « in Crosland's opinion, the categories into which scientists are divided by prosopographical studies often include many which have no significance for the purely scientific work of their subjects […] and that prosopography runs the danger of treating its subjects as if they were of equal importance » ; lui-même insiste sur le fait que si son étude met les membres de l'Académie des sciences au centre de l'attention, ce n'est pas « simply as the objects of prosopography[4] ». Ce qui commence par un malaise — comment traiter « comme s'ils étaient d'égale importance » les créateurs célèbres et les presque inconnus — transforme parfois l'enquête prosopographique en un simple outil de contraste, pour étudier à travers la moyenne d'un groupe social la nature exceptionnelle de certains en tant qu'ils échappent à cette moyenne, se *distinguent*. Évoquant ce type d'usages en histoire de la technologie, Christine Mac Leod et Alessandro Nuvolari soulignent pourtant qu'il contribue paradoxalement à surévaluer le rôle des inventeurs individuels et à maintenir en fin de compte une vision héroïque du changement et de l'innovation[5].

---

[1] Koenraad VERBOVEN, Myriam CARLIER, Jan DUMOLYN, « A Short Manual to the Art of Prosopography », in Katherine S. B. KEATS-ROHAN (ed.), *Prosopography Approaches and Applications: A Handbook*, Oxford : Occasional Publication of the Unit for Prosopographical Research, Linacre College, University of Oxford, 2007, p. 35-69 ; version en ligne, https://biblio.ugent.be/record/376535.

[2] Voir la définition donnée dans l'article classique de Lawrence STONE, « Prosopography », *Daedalus* (100), 1971, p. 46-79 : « Prosopography is the investigation of the common background characteristics of a group of actors in history by means of a collective study of their lives. The method employed is to establish a universe to be studied, and then to ask a set of uniform questions. The various types of information about the individuals in the universe are then juxtaposed and combined and are examined for significant variables. »

Je voudrais proposer ici un usage différent des prosopographies : il s'agit de les utiliser d'emblée pour appréhender un seul individu, qui plus est, un mathématicien renommé, dont les résultats ont reçu les louanges de ses contemporains tout autant que des mathématiciens actuels. De telles prosopographies me semblent en fait un outil puissant pour éviter deux des écueils ordinaires de la biographie évoqués plus haut : d'une part, la question des sources pertinentes et de leur traitement, d'autre part, le piège de la focalisation individuelle, qui tend à capter au seul profit de la figure centrale toutes les pratiques de la création détectées.

L'exemple choisi ici est celui de Charles Hermite. « Parmi les mathématiciens de tous les temps, il en est peu qui aient exercé une influence directe comparable à celle d'Hermite ; il n'en est pas dont l'œuvre soit plus sûrement impérissable[6] », écrit à sa mort Paul Painlevé, lui-même membre de l'Institut, futur ministre et Président du Conseil. Né à Dieuze, en Lorraine, en 1822, Hermite a été élu à l'Académie des sciences en juillet 1856, maître de conférences à l'École normale supérieure en 1862, puis professeur à l'École polytechnique et à la Faculté des sciences de Paris dans les années 1870. Au centre des institutions françaises, il a entretenu d'intensives correspondances internationales, de la Suède à l'Italie, de la Russie aux États-Unis. Cette gloire n'est pas éphémère : depuis l'année 2000, plus de mille articles recensés dans *Math Sci Net* contiennent le mot « Hermite » dans leur seul titre. Il existe (et cette liste n'est pas exhaustive) des polynômes et une constante d'Hermite, des opérateurs et des matrices d'Hermite, des approximants de Hermite-Padé, des inégalités de Hermite-Hadamard, sans oublier les formes hermitiennes.

Plusieurs essais biographiques ont naturellement été consacrés à Hermite[7]. Ce que je propose ici, c'est de l'aborder par le biais de deux enquêtes prosopographiques. Mais il ne s'agit nullement de *le* mettre en série avec d'autres, que ce soit pour en récupérer les traits moyens ou pour le distinguer. Au contraire, les sériations, multiples, seront ici mises à profit pour comprendre des phénomènes directement biographiques, pour capturer des pratiques en tant qu'elles sont propres à Hermite lui-même.

---

3       Pierre BOURDIEU, *Les Règles de l'art*, Paris : Seuil, 1992, p. 10.
4       David J. STURDY, *Science and Social Status. The Members of the Académie des sciences, 1666-1750*, Woodbridge : Boydell Press, 1995, p. xiv.
5       Christine MAC LEOD et Alessandro NUVOLARI, « The Pitfalls of Prosopography: Inventors in the Dictionary of National Biography », *Technology and Culture* (47), 2006, p. 757-776.
6       Paul PAINLEVÉ, « Ch. Hermite », *La Nature* (1445), 2 février 1901, p. 145-146.
7       Voir en particulier Elena Petrovna OZHIGOVA, *Charles Hermite : 1822-1901* (en russe), Leningrad: Nauka, 1982; Claude BREZINSKI, *Charles Hermite. Père de l'analyse mathématique moderne*, Paris: SFHST, 1990.

I-Une prosoprographie d'articles

La première enquête dont nous parlerons est en fait double : elle s'appuie sur une prosopographie d'articles mathématiques *et* sur une prosopographie de personnes. Les articles sont les quelque 500 consacrés à la théorie des nombres et parus entre 1870 et 1914 dans des journaux publiés en France ; les personnes sont toutes celles mentionnées dans ces articles, principalement leurs auteurs ou ceux cités en référence à l'intérieur des articles.

L'objectif initial de cette enquête ne concernait pas Hermite, mais bien l'historiographie de la théorie des nombres[8]. Le statut de ce domaine change radicalement au cours du $19^e$ siècle : marginal, relativement élémentaire et quelque peu disparate au tout début du siècle, il devient à la fin une branche majeure des mathématiques, moins d'ailleurs par son importance numérique réelle(seulement 3-4% des articles de mathématiques publiés relèvent de ce domaine vers 1900) que par ses liens à de multiples autres branches (analyse complexe, algèbre, géométrie, théorie des ensembles, …) et par son rôle de modèle théorique. Un tiers des problèmes que David Hilbert propose à l'attention des mathématiciens du $20^e$ siècle au Congrès international des mathématiciens de 1900 relève de ce domaine. Or, dans les années 1990 encore, l'histoire de la théorie des nombres était une histoire perçue comme presque exclusivement allemande, focalisée sur l'apport successif de quelques auteurs célèbres, Carl Friedrich Gauss, Carl Gustav Jakob Jacobi, Peter Gustav Dirichlet, Ernst Eduard Kummer, Richard Dedekind, David Hilbert. Il s'agissait donc alors pour moi de comprendre comment, dans ces conditions, la production d'articles de théorie des nombres dans les journaux français pouvait être si abondante, et pour cela, d'examiner systématiquement l'ensemble de cette production afin de mieux comprendre les phénomènes en jeu.

Comme pour les prosopographies usuelles, tous ces articles ont été intégrés dans une base de données[9]. Y figurent, outre les données bibliographiques de chaque article (permettant en particulier des recoupements par journal), les sections du classement de l'article dans plusieurs taxonomies mathématiques en usage à l'époque considérée, l'ensemble des références utilisées dans l'article (repérées par leurs auteurs et les lieux de publication), et des mots-clés indexant en particulier les méthodes, les résultats, certains concepts.

Inscrit dans le genre bien balisé des études quantitatives sur les sciences, ce travail pose d'ailleurs des problèmes de mise en œuvre et d'interprétation à peu près analogues à ceux des prosopographies de personnes[10], en premier lieu la sélection même des entrées de la prosopographie, c'est-à-dire ici des articles.

À partir de 1868, les mathématiciens disposent d'un recensement quasi-systématique de tous les articles parus, classés par sujets, le *Jahrbuch über die Fortschritte der Mathematik* ; la théorie des nombres, « Zahlentheorie », est une des rubriques proposées, ce qui semble fournir une représentation objective et synchrone de ce domaine, induisant un choix naturel des articles retenus. Or, on constate que des articles jugés maintenant essentiels dans le développement de la théorie des nombres n'apparaissent pas sous cette rubrique dans le *Jahrbuch* ; ils sont classés par exemple dans des rubriques de géométrie ou d'analyse. Mais il y a plus : d'autres sources ou classements contemporains fournissent également des listes d'articles différents. Dans les trois volumes de *The History of the Theory of Numbers*, compilés par Leonard Dickson et ses collaborateurs dans les premières décennies du $20^e$ siècle, par exemple, sont inclus des articles pourtant non classés en théorie des nombres dans le *Jahrbuch* ; le *Répertoire bibliographique*, outil concurrent au *Jahrbuch* proposé à la fin du $19^e$ siècle, adopte un classement encore différent ; les journaux eux-mêmes indexent dans certains cas leurs articles, et là encore leurs classifications ne sont pas les mêmes.

---

8  Elle est présentée dans Catherine GOLDSTEIN, « La théorie des nombres dans les Notes aux *Comptes rendus de l'Académie des sciences* (1870–1914) : un premier examen », *Rivista di Storia della scienza* (II, 2, 2), 1994, p. 137-160 et « Sur la question des méthodes quantitatives en histoire des mathématiques : le cas de la théorie des nombres en France (1870--1914) », *Acta historiae rerum naturalium nec non technicarum* (New series 3), 1999, p. 187-214.

9  Élaborée de manière artisanale à l'aide du logiciel Hypercard, cf. Catherine GOLDSTEIN, « Sur la question des méthodes... », article cité note 8, elle est en cours de transfert dans la base THAMOUS initiée par Alain Herreman.

10  Ces problèmes sont examinés en détail dans Catherine GOLDSTEIN, « Sur la question des méthodes... », article cité note 8.

Même en retenant tous les articles classés en théorie des nombres par l'un de ces classements, le fait de les mettre sur le même plan peut soulever des objections. Les modes de publication dans les différents journaux (*Journal de mathématiques pures et appliquées*, *Comptes rendus hebdomadaires des séances de l'Académie des sciences*, *Annales de l'École normale supérieure*, *Journal de l'École polytechnique*, *Bulletin de la Société mathématique de France*, *Nouvelles Annales de mathématiques*, pour ne citer que les plus significatifs dans ce domaine) sont en effet extrêmement variés. Les notes aux *Comptes rendus* sont limitées en nombre de pages : au début de la période considérée, certains auteurs découpent leur travail en petites unités, les publiant à la file dans les *Comptes rendus* ; plus tard, au contraire, la norme est plutôt de publier un court résumé des résultats obtenus dans les *Comptes rendus*, alors que les détails et les preuves sont données dans une autre revue. Ces comportements se chevauchent parfois dans le temps, voire apparaissent en alternance chez un même auteur. Faut-il donc, et dans quel cas, inclure chaque note séparément ou non ? Par ailleurs, les lecteurs visés par les *Nouvelles Annales* sont les étudiants préparant les grandes écoles ou leurs enseignants : certains articles de ce journal sont effectivement des articles de recherche originaux, même si les méthodes sont élémentaires, mais d'autres sont exclusivement destinés à promouvoir un exercice ou une technique dans l'enseignement, sans prétention à l'originalité, d'autres encore visent à diffuser une approche développée à l'étranger. Que faut-il prendre en compte ? Le *Journal de mathématiques pures et appliquées* ou les *Annales de l'École normale supérieure* eux-mêmes ne publient pas que des articles de recherche nouveaux : certaines traductions de journaux de recherche étrangers, ou des lettres y sont insérées. Fait-il sens de comparer ces différents types de contributions par un ensemble d'indicateurs communs ?

Ces questions se posent bien sûr, de manière plus attendue, pour l'autre prosopographie de cette enquête, celle de la base des personnes. Y sont indiquées la fonction, ou les fonctions, de ces personnes *dans la base d'articles* (selon qu'ils y apparaissent comme auteur d'articles de la base, auteur d'une référence citée, traducteur, etc.), ainsi que des informations biographiques et institutionnelles normalisées. Même si on se restreint à une seule fonction, comme celle d'auteur d'articles de la base, se côtoient ainsi des célébrités du domaine, comme Richard Dedekind lui-même (qui publie en français, dans le *Journal de mathématiques pures et appliquées*, un article fondateur pour l'avènement d'un point de vue structural en théorie des nombres) ou Henri Poincaré (qui élabore les liens entre formes quadratiques ternaires et groupes fuchsiens entre le *Journal de l'École polytechnique* et les *Comptes rendus hebdomadaires des séances de l'Académie des sciences*), et des amateurs, enseignants comme Honoré Desboves ou militaires comme Ernest Laquière, qui, même s'ils s'intéressent en général à des questions techniquement élémentaires, peuvent être eux aussi néanmoins extrêmement productifs[11].

Ces deux bases permettent d'obtenir des informations sur la théorie des nombres publiée en France à la fin du 19$^e$ siècle d'au moins deux façons : soit de manière statistique, directement à partir des données, disjointes mais uniformisées des bases, soit en établissant des liens entre données à l'intérieur de ces bases, afin de structurer les informations[12]. En considérant par exemple comme un lien la référence dans un article à un autre, ou encore l'usage d'un même mot-clé, j'ai pu construire des réseaux (*clusters*) de textes, c'est-à-dire des sous-populations de la base de textes : à l'intérieur d'un même réseau, les articles font référence les uns aux autres, ils partagent aussi un noyau de sources communes. En revanche, deux réseaux différents sont à peu près étanches du point de vue des citations : même lorsqu'ils semblent traiter de la même question, deux articles de réseaux distincts ne font pas allusion l'un à l'autre. Lorsqu'il y a occasionnellement référence à un autre réseau, c'est le plus souvent pour en critiquer les résultats ou l'approche, ou peut-être revendiquer une priorité mal reconnue. La structure de ces réseaux est complexe, avec des sous-réseaux partiels, et quelques rares fusions au cours de la période observée.

---

11    Voir sur Desboves par exemple, Norbert SCHAPPACHER, « Le développement de la loi de groupe sur une cubique », in *Séminaire de théorie des nombres de Paris 1988-1898*, Progress in Mathematics 91, Basel, Boston: Birkhäuser, 1991, p. 159-184.

12    Voir par exemple sur la première approche, Hélène GISPERT, *La France mathématique : la Société mathématique de France (1870-1914)*, Paris : SFHST et SMF, 1991, en particulier chap. 7, et sur la seconde, la construction de réseaux de textes expliquée dans les articles cités note 8.

De façon analogue, la base de personnes s'organise en réseaux variés, qu'ils soient définis à partir de liens personnels (existence de correspondance, traces de visites, de cours reçus ou donnés), ou de ceux-là même définis par la base de textes (le lien entre auteurs y étant défini par le fait que l'un cite l'autre dans un article de la base). Il est remarquable que les réseaux personnels et réseaux de textes ne coïncident pas nécessairement (ce qui justifie d'ailleurs *a posteriori* l'intérêt d'avoir recours à ces deux bases prosopographiques). Pour n'évoquer qu'un exemple, James Joseph Sylvester est un ami personnel de Charles Hermite—ils échangent des lettres, Hermite est intervenu auprès de collègues étrangers pour obtenir des lettres de recommandation lorsque Sylvester cherchait un poste— et dans les années 1850, ils ont tous deux contribué, avec Arthur Cayley, à l'élaboration d'un nouveau champ de recherches mathématiques, la théorie des invariants. Sylvester est aussi un auteur dans la base de textes (sur la période 1870-1914) : or, il ne cite pas les articles d'Hermite de cette base, et qui plus est ses textes ne font pas partie du réseau faisant référence aux articles de théorie des nombres d'Hermite (antérieurs à la période prise en compte dans la base).

Néanmoins, même si les réseaux de textes diffèrent de ceux obtenus pour leurs auteurs à partir des liens personnels, ils peuvent présenter des propriétés *biographiques,* c'est-à-dire mettent en évidence des caractèresrelatifs à leurs auteurs au-delà du simple fait de faire référence les uns aux autres. Certains réseaux d'articles par exemple font apparaître une participation à une société savante, comme la Société mathématique de France ou l'Association française pour l'avancement des sciences. D'autres témoignent de solidarités internationales particulières : certains réseaux citent surtout des articles parus dans des revues américaines ou françaises, d'autres des articles parus dans des revues allemandes (un trait significatif après la guerre franco-prussienne). Pour certains réseaux, tous les auteurs sont professeurs d'université, pour d'autres, les professions sont plus variées et non significatives, etc. Une autre propriété intéressante est celle des modes de publication : certains réseaux d'articles sont caractérisés par des rafales de petits articles ou de notes, d'autres sont constitués surtout de longs articles parus dans un journal particulier.

Par ailleurs, certains auteurs ne publient que des articles appartenant à un seul réseau de textes, alors que d'autres publient dans plusieurs réseaux de textes disjoints. Ce dernier point montre d'ailleurs de manière frappante la distinction entre les réseaux obtenus à partir des deux bases de données, réseaux de personnes et réseaux de textes : l'étanchéité entre deux réseaux d'articles disjoints ne cesse pas alors même que certaines personnes sont auteurs d'articles dans les deux. Ceci veut dire, *par définition du réseau d'articles*, que le comportement de l'auteur est tel que l'article qu'il publie ne se distingue en rien des autres articles du même réseau : il cite les mêmes textes que les autres articles du réseau, et ne cite pas ceux de l'autre réseau (y compris les siens). Cette situation est illustrée par Edmond Maillet, ingénieur des Ponts et Chaussées, prix Poncelet de l'Académie des sciences en 1912, qui écrit sur de multiples sujets arithmétiques à la fin du 19e siècle[13].

Revenons maintenant à Charles Hermite. Celui-ci était considéré dans l'historiographie comme un des rares mathématiciens non allemands à avoir contribué de manière significative au développement de la théorie des nombres. André Weil écrit ainsi de manière caractéristique : « The great number-theorists of the last century are a small and select group of men. The names of Gauss, Jacobi, Dirichlet, Kummer, Hermite, Eisenstein, Kronecker, Dedekind, Minkowski, Hilbert spring to mind at once[14]. » Hermite est l'auteur de résultats célèbres sur les valeurs entières de formes (inventant au passage la notion des formes dites maintenant hermitiennes), sur l'approximation par des nombres rationnels, et, en 1873, du premier théorème démontrant la transcendance d'une constante usuelle de l'analyse, *e* en l'occurrence[15]. Or, malgré cette importance reconnue, Hermite n'apparaît que d'une manière anecdotique parmi les auteurs d'articles de notre première base : sont seulement inscrits dans la base un court article de lui dans les *Nouvelles Annales* en 1872 sur les solutions rationnelles de l'équation $x^3+y^3=z^3+v^3$ et deux publications de quelques pages sur les formes quadratiques et les fractions continues dans le *Bulletin des sciences mathématiques* dans les années 1880.

---

13    Les auteurs prolifiques à l'intérieur d'un *même* réseau d'articles semblent d'ailleurs plus susceptibles de faire l'objet d'une biographie : par exemple, Édouard Lucas (Anne-Marie DECAILLOT, *Édouard Lucas (1842–1891), Le parcours original d'un scientifique français dans la deuxième moitié du 19ᵉ siècle*, Thèse de l'université Paris V, 1999) ; Jacques Hadamard (Vladimir MAZ'YA, Tatiana SHAPOSHNIKOVA, *Jacques Hadamard, A Universal Mathematician*, History of Mathematics 14, Providence: AMS and London : LMS, 1998) ; Albert Châtelet (Jean-François CONDETTE, *Albert Châtelet, La République par l'école (1883-1960)*, Arras: Artois Presses Université, 2009), … On voit ainsi réapparaître quelques-uns des critères traditionnels séparant les sujets qui relèvent du genre biographique de ceux qui relèvent de la prosopographie.

Cette simple absence n'est pas interprétable immédiatement à partir de la prosopographie donnée, mais elle a une valeur heuristique. Elle incite à l'examen séparé des publications d'Hermite sur cette période, examen qui fait apparaître plusieurs phénomènes distincts. D'abord, les recherches d'Hermite s'orientent après les années 1870 vers des sujets d'analyse, en particulier les équations aux dérivées partielles ; ce sont des sujets qu'il a toujours étudiés, mais dont les liens avec l'algèbre et la théorie des nombres, essentiels dans les années 1850-1860, sont ensuite moins apparents. Ensuite, Hermite publie ses travaux de théorie des nombres surtout hors de France (dans le *Journal für die reine und angewandte Mathematik,* les *Annali di matematica pura ed applicata, Acta Mathematica,* entre autres*)*. Enfin — et ceci confirme les difficultés déjà notées pour établir notre prosopographie de textes —, certains de ses articles importants pour la théorie des nombres, dont ses notes sur la transcendance de *e,* sont classées en analyse («BesondereFunkt.*»* par exemple pour les notes sur *e*) et non en « Zahlentheorie » dans le *Jahrbuch*. C'est la conjonction de ces trois effets qui explique le faible nombre des articles d'Hermite dans la base prosopographique de textes.

Mais qu'en est-il d'Hermite dans la base des personnes ? Rappelons que n'y sont pas relevés seulement les auteurs des articles de la base prosopographique de textes, mais tous les noms propres qui y interviennent (par exemple les auteurs de travaux cités dans ces textes). Une recherche globale du nom « Hermite » dans la base montre qu'il a occupé d'autres rôles au sein de la théorie des nombres entre 1870 et la première guerre mondiale que celui d'auteur de quelques textes.

Il intervient d'abord comme une des principales *références* communes à un important réseau de textes. Plus précisément, ce réseau est caractérisé par le fait que les articles le composant se réfèrent, outre leurs citations réciproques, à la section V des *Disquisitiones arithmeticae* de Gauss et à des travaux d'Hermite ou de Leopold Kronecker sur les formes[16]. Entre 1870 et 1914, ce réseau de textes réunit environ 90 contributions, publiés surtout dans les *Comptes rendus hebdomadaires des séances de l'Académie des sciences* et le *Journal de mathématiques pures et appliquées*. Un thème partagé dans ce réseau d'articles est l'étude des formes *m*-aires de degré *n* (c'est-à-dire des polynômes homogènes de degré *n* à *m* variables), soit à coefficients entiers, soit dont on considère les valeurs entières, les classifications de ces formes, leurs invariants, leur réduction (par transformations linéaires) à des formes plus simples à petits coefficients, leurs liens avec certaines fonctions complexes ou avec des problèmes d'approximation. Les auteurs les plus actifs d'articles de ce réseau[17] sont Henri Poincaré, Camille Jordan, Georges Humbert, Jean-Armand de Séguier,

---

[14] André WEIL, « Introduction », in *Collected Papers. Ernst Eduard Kummer*, New York: Springer, 1975, citation p. 1. On remarque, comme indiqué plus haut, que l'histoire du domaine est réduite à celle de l'œuvre de quelques figures importantes. On trouvera un point de vue différent dans Catherine GOLDSTEIN, Norbert SCHAPPACHER & Joachim SCHWERMER, eds., *The Shaping of Arithmetic after C. F. Gauss's* Disquisitiones arithmeticae, New York, Berlin : Springer, 2007.

[15] Sur les travaux d'Hermite sur les formes, voir Catherine GOLDSTEIN, « The Hermitian Form of Reading the *Disquisitiones* », in *The Shaping of Arithmetic, op. cit.* note 14, p. 377-410. Sur la preuve de la transcendance de *e*, voir Michel WALDSCHMIDT, « Les débuts de la théorie des nombres transcendants », *Cahier du Séminaire d'histoire des mathématiques* (4), 1983, p. 93-115, ainsi que la présentation des articles d'Hermite par Michel Waldschmidt sur le site Bibnum, http://www.bibnum.education.fr/mathématiques/théorie-des-nombres/la-démonstration-de-la-transcendance-de-e#.

[16] Ce réseau est pour cette raison appelé « cluster H-K» dans *The Shaping of Arithmetic, op. cit.* note 14.

[17] En étendant ces bases à la production étrangère, on pourrait inclure des auteurs comme Eduard Selling, Paul Bachmann, Hermann Minkowski, Henry John Stephen Smith, Luigi Bianchi, etc.

Léon Charve, Albert Châtelet, …. Le nom d'Hermite apparaît en référence dans leurs textes à cause de ses résultats des années 1850-1865 sur les formes, les transformations linéaires associées et les fractions continues. Châtelet mentionne par exemple[18] : « Hermite a énoncé cette propriété pour une substitution linéaire d'ordre quelconque, *Journal de Crelle* 41 ». Charve rappelle en 1880 que « Dans ses Lettres à Jacobi, où sont exposés tant de beaux résultats, M. Hermite est revenu plusieurs fois sur la théorie de la réduction[19] ». Ou encore, dans un article de Poincaré de 1882 : « Je résoudrai ces problèmes par une généralisation de la méthode de M. Hermite, sur laquelle je veux donner d'abord quelques explications[20] ». Ce réseau mesure donc en partie les traces de son impact scientifique en théorie des nombres à la fin du 19e siècle.

Hermite apparaît également dans une autre fonction : c'est lui qui présente presque toutes les notes de théorie des nombres insérées dans les *Comptes rendus des séances de l'Académie des sciences*. Son nom apparaît donc au début de la note (voir Illustration 1). Les notes aux *Comptes rendus* jouent un rôle international grandissant dans la deuxième moitié du siècle, de nombreux mathématiciens étrangers profitant de leur diffusion rapide et large pour annoncer leurs travaux par leur intermédiaire à l'ensemble de la communauté mathématique. À ce titre, Hermite se trouve donc au centre d'une structure de diffusion essentielle du domaine, indépendante des réseaux d'articles définis précédemment dans lesquels telle ou telle note peut se trouver. Il présente ainsi à l'Académie des notes de Jacques Hadamard sur la fonction zeta de Riemann ou de Théophile Pépin sur des équations diophantiennes (Illustration 1), articles de notre base qui ne le citent pas et ne font pas partie du réseau Hermite-Kronecker (cluster H-K, voir note 16). Il s'agit donc d'un aspect de son impact professionnel, au-delà d'une influence mathématique directe marquée par la référence à ses recherches.

ANALYSE MATHÉMATIQUE. — *Nouveaux théorèmes sur l'équation indéterminée* $ax^4 + by^4 = z^2$. Note du P. PEPIN, présentée par M. Hermite.

« Les nouveaux théorèmes renfermés dans cette Note sont semblables à ceux que j'ai eu l'honneur de présenter à l'Académie sur le même sujet ;

Illustration 1 : Le rôle d'Hermite pour les notes de théorie des nombres dans les *Comptes rendus* de l'Académie des sciences. Cet article est classé par le *Jahrbuch* en « Zahlentheorie » ; on remarque le classement différent des *Comptes rendus* (« analyse mathématique »).

Ces deux fonctions n'épuisent pas les apparitions du nom « Hermite » dans ces bases prosopographiques. Plusieurs articles de la base sont des extraits de lettres dont il est le destinataire[21], il est aussi mentionné pour avoir encouragé la traduction d'un article, ou simplement un travail particulier ; Poincaré indique par exemple dans une note sur la réduction des formes : « D'après les conseils de M. Hermite, j'ai poursuivi les résultats obtenus et j'ai cherché à approfondir les conditions d'équivalence ou des substitutions semblables de pareils systèmes[22] ».

Autrement dit, à partir de ces sources sérielles, au départ des éléments d'articles mathématiques (pour la plupart d'autres auteurs qu'Hermite), il nous a été possible de capturer des aspects singuliers, propres à lui : Hermite est apparu comme père fondateur d'une communauté de recherches (définie par un réseau d'articles lié par des références réciproques), il est aussi apparu comme un intermédiaire privilégié, plaque tournante internationale pour la publication.

---

18  Albert CHÂTELET, « Contribution à la théorie des fractions continues arithmétiques», *Bulletin de la Société mathématique de France* (40), 1912, p. 1-25, citation p. 11.
19  Léon CHARVE, « De la réduction des formes quadratiques ternaires positives et de son application aux irrationnelles du troisième degré », *Annales scientifiques de l'École Normale Supérieure* (Sér. 2, 9), 1880, p. 3-156, citation p. 8.
20  Henri POINCARÉ, « Sur les formes cubiques ternaires et quaternaires », *Journal de l'École polytechnique* (51), 1882, p. 45-91, citation p. 45.
21  Par exemple Rudolph LIPSCHITZ, « Propositions arithmétiques tirées de la théorie de la fonction exponentielle. (Extrait d'une Lettre adressée à M. Hermite) », *Journal de Mathématiques Pures et appliquées* (4e s. 2), 1886, p. 219-238.
22  Henri POINCARÉ, « Sur la réduction simultanée d'une forme quadratique et d'une forme linéaire », *Comptes rendus de l'Académie des sciences* (91) , 22 novembre 1880, p. 844-846.

Il va de soi que nous aurions aussi pu appréhender Hermite en tant qu'« homme moyen », contribuant aux effets collectifs attendus d'une prosopographie. Par exemple, j'ai indiqué plus haut que les auteurs d'articles de certains réseaux ont des métiers identiques, par exemple des universitaires ; grâce à ses quelques articles (qui relèvent d'un de ces derniers réseaux), Hermite, professeur d'université dans les années 1880, contribue donc à la *norme* du réseau d'article dont il est un auteur. Il est aussi possible qu'en prenant en compte la production des années 1840-1870 en théorie des nombres (ou bien un autre domaine des mathématiques), on puisse récupérer, par contraste cette fois avec la moyenne, certains caractères qui rendraient Hermite *banalementexceptionnel*, peut-être par exemple un nombre d'articles de théorie des nombres élevé par rapport aux autres auteurs français.

Mais l'usage des prosopographies proposé ici est différent. Il ne s'agit ni de noyer Hermite au milieu d'autres mathématiciens pour lesquels on aurait recueilli les mêmes informations, ni de l'en faire ressortir par son écart à la norme. L'usage est transverse — un type d'usage facilité par l'informatisation des ressources, qui fournit des fonctions « recherche » diversifiées et efficaces. Il donne accès à des caractères exceptionnels (dans le sens de « spécifiques ») d'un scientifique, à ses positions particulières au sein d'une collectivité. Le même genre de transversalité est parfois mis en l'œuvre pour dégager des sous-communautés dans la communauté de référence qui fait l'objet de la prosopographie : ici, « Hermite » est bien traité techniquement comme une sous-communauté. Cette possibilité me semble soulever deux points importants : d'abord que de multiples analyses prosopographiques peuvent être pertinentes, que l'on s'intéresse à une personne ou même à un collectif de personnes ; ensuite, qu'il n'y a pas de sources biographiques naturelles. Nous avons eu recours ici à des articles écrits par d'autres personnes, pour y saisir ensuite le nom d'Hermite. Les biographies jouent souvent sur un effet de réalisme induit par des sources dont l'adéquation au sujet est comme postulée à l'avance. Un détachement préliminaire entre sources et personne étudiée peut permettre en particulier d'aborder plus efficacement le rapport de la vie à l'œuvre et de pallier en partie le défaut de perspective évoqué au début de cet article.

II— Des noms propres dans les *Œuvres* de Charles Hermite

Prenant le contre-pied de la première, ma deuxième enquête prosopographique est entreprise cette fois à partir de sources qu'il est plus banal d'attacher à Hermite : ses *Œuvres*, dans la version posthume publiée en quatre volumes « sous les auspices de l'Académie des sciences » entre 1905 et 1917, sous la direction d'Émile Picard (mathématicien, académicien et gendre d'Hermite)[23]. Je ne m'intéresserai pas ici pourtant aux résultats mathématiques qu'elles contiennent, mais aux personnes qui y sont mentionnées ; en ce sens, cette deuxième enquête s'inscrit dans la continuité de la première. J'ai exclu les nécrologies et d'autres textes de circonstance qui abondent en noms propres (comme le discours d'Hermite lors de l'inauguration de la Nouvelle Sorbonne le 5 août 1889), car il s'agissait dans un premier temps de se concentrer sur les noms propres volontairement mobilisés par Hermite dans sa pratique scientifique. On remarquera que cette restriction limite d'emblée les contours sociaux d'un Hermite scientifique, face à un éventuel Hermite institutionnel : les problèmes de sources se posent donc tout autant pour cette approche que dans une approche biographique plus traditionnelle, même s'ils ne sont pas les mêmes.

Restent à prendre en compte environ 190 écrits, articles de journaux pour la plupart et quelques sections de livres. J'y ai relevé systématiquement tous les noms propres, qu'ils apparaissant dans le texte ou dans les paratextes (titres ou notes par exemple, à l'exception des notes plus tardives ajoutées par les éditeurs des *Œuvres*). Une base de données les rassemble sous la forme sous laquelle ils apparaissent (avec ou sans avant-noms, avec ou sans qualificatifs, …), et inclut le lieu de publication, l'année, les mentions qui entourent leurs occurrences (voir Illustration 2). À l'intérieur d'un même article, chaque nom distinct a été relevé une seule fois (en précisant les différents contextes d'apparition si nécessaire) ; un autre choix semblait peu pertinent, car, à l'intérieur d'un même article, Hermite reprend souvent par des périphrases des noms de personnes déjà mentionnées. Par exemple, dans une série de notes[24] aux *Comptes rendus* en 1855, après avoir évoqué les résultats d'Adolf Göpel et de Johann Georg Rosenhain sur les fonctions hyperelliptiques, Hermite précise que « ces illustres géomètres ont en même temps donné … l'expression analytique de treize autres fonctions… » ; une expression analogue (« illustre géomètre ») sert ensuite dans les mêmes notes à désigner Göpel seul, puis Abel. Compter les différentes mentions d'un nom dans un même article n'était donc évident ni à définir, ni à interpréter, et j'ai choisi de négliger cette question ici. Cette convention de relevé fournit un peu plus de 800 occurrences, réparties entre 159 noms propres distincts[25] ; une base secondaire inclut ensuite des informations extérieures sur ces 159 personnes (prénom, dates de vie, lieux d'exercices professionnels, thèmes de recherches, …).

---

23  Charles HERMITE, *Œuvres*, éd. Émile Picard, 4 volumes, Paris : Gauthier-Villars, 1905-1917.

Parmi eux, 15 sont morts avant la fin du 18$^e$ siècle : Francis Bacon, Jacob Bernoulli, Jean Le Rond d'Alembert, René Descartes, Euclide, Leonhard Euler, Pierre Fermat, Johannes Kepler, Johann Lambert, John Landen, Colin MacLaurin, Isaac Newton, Michel Rolle, James Stirling, Brook Taylor. À part le premier, mentionné par Hermite pour avoir dit que l'admiration est le principe du savoir, ces noms n'interviennent la plupart du temps que pour désigner des objets et des résultats mathématiques[26] : « équation de Kepler » (*Proceedings of the London Mathematical Society* VII, 1876), « formule de MacLaurin » (*Annali di Matematica pura ed applicata*, 2$^e$, III, 1869), « série de Stirling » (*American Journal of Mathematics* 17, 1895), « nombres de Bernoulli » (*Journal für die reine und angewandte Mathematik* 81, 1876), etc. La fonction adjectivale de ces noms est d'ailleurs soulignée par la coexistence d'expressions comme « intégrale eulérienne » ou « processo Euclideo ». Notons enfin que cet emploi n'est pas limité aux auteurs anciens, même s'il n'est alors plus dominant : on trouve ainsi des « transcendantes de Bessel », des « séries de Fourier », des « polynômes de Legendre », une « intégrale de Poisson ».

Douze auteurs du 19$^e$ siècle sont utilisés comme compléments du nom « journal » ou de variantes : « archives de Grunert », « bulletin de M. Darboux », « annales de M. Tortolini ». Cet usage peut sembler très proche du précédent, correspondant alors à une totale désincarnation du nom. Mais deux indices montrent que ce n'est pas le cas. Un indice linguistique, d'abord, la présence de « M. » en avant-nom dans les titres des journaux[27]. En général, dans tout le corps des textes, cette présence distingue simplement les vivants des morts : « M. Cauchy », dans les *Comptes rendus des séances de l'Académie* en 1851, ou dans le *Journal für die reine und angewandte Mathematik* en 1854, devient ainsi « Cauchy » tout court dans les *Annali di matematica pura ed applicata*, 2$^e$ s., III, en 1869 ou le *Bulletin des sciences mathématiques* en 1879 (Cauchy est mort en 1857)—de rares exceptions ont lieu dans le cadre d'hommages particuliers, réactivant la personne même. Or, on trouve de manière analogue un « journal de M. Thompson » (autrement dit *Cambridge and Dublin Mathematical Journal*) en 1854, mais des « Annales de Gergonne » (sans M.) lorsqu'Hermite fait référence dans les années 1880 à ce journal, autrement dit les *Annales de*

---

24 Charles HERMITE, « Sur la théorie de la transformation des fonctions abéliennes », *CRAS* (XL), 1855, rep. in Charles HERMITE, *op. cit.* n. 23, vol. I, p. 444-478, citation p. 444.
25 En dehors du nom d'Hermite lui-même. Remarquons toutefois que certaines indications biographiques accompagnent son nom lors de ses premiers articles : « M. Hermite, élève au collège Louis-le-Grand (institution Mayer) » pour ses premières contributions aux *Nouvelles Annales*, en 1842 ; « M. Hermite, élève de l'École polytechnique », dans les *Mémoires présentés par divers savants à l'Académie des sciences*, t. X.
26 Les rares exceptions sont quelques occurrences d'auteurs du 18$^e$ siècle : (Jacob) Bernoulli est cité 2 fois en 1878, certaines formules sont attribuées à Euler en 1872 (*Nouvelles Annales* 2$^e$ s., 11) et Lambert, qui apparaît dans l'expression « fraction continue de Lambert », intervient aussi (2 occurrences) comme auteur historique de l'unique preuve alors connue de l'irrationalité de $\pi$ et $\pi^2$ dans des articles de 1873 où Hermite démontre la transcendance de *e*.
27 Il reste difficile de contraster cette situation avec l'usage des noms dans la désignation d'objets mathématiques, car peu portent le nom d'auteurs vivants dans l'œuvre d'Hermite. On trouve quand même une « réduite de Jerrard » dans les années 1850 (Charles HERMITE, *op. cit.* n. 23, vol. 2, p. 82), conjointement à des « M. Jerrard » (Georg B. Jerrard meurt en 1863).

*mathématiques pures et appliquées* ; Joseph Gergonne est mort en 1859, William Thomson (Lord Kelvin) ne décèdera qu'en 1907.Un deuxième indice est le soin apporté à suivre les changements d'éditeurs : dans les articles d'Hermite, le « journal de M. Crelle » fait place au « journal de M. Borchardt », puis au « journal de MM. Weierstrass et Kronecker » après 1881, alors qu'il s'agit dans tous les cas de ce que nous appelons maintenant le *Journal für die reine und angewandte Mathematik*[28]. De même pour le « journal de M. Liouville » devenu « journal de M. Jordan ». À l'intérieur d'un même article (par exemple en 1888, dans le *Bulletin de l'académie des sciences de St Pétersbourg*), deux noms propres associés au même journal peuvent d'ailleurs cohabiter, par exemple si les références sont faites à des articles d'années différentes, donc parus dans des volumes édités par des personnes différentes (en l'occurrence Crelle et Borchardt[29]).

Ces deux indices suggèrent une forte présence personnelle des éditeurs de journaux. L'usage de leur nom au lieu du titre officiel du journal (ou dans certains cas, en sus du titre) est donc significatif ici. Ceci conforte bien sûr la perception des historiens des journaux mathématiques qui insistent tous sur le rôle central et structurant joué par les éditeurs successifs[30]. Mais vu d'un mathématicien qui ne fut pas lui-même éditeur principal d'une revue, ce phénomène prend une connotation plus intime : une simple référence semble évoquer non seulement son auteur, mais aussi d'autres acteurs du travail mathématique, les éditeurs des journaux où cette référence a été publiée, et ce d'une manière tout aussi individualisée. Cette conclusion est corroborée par un autre fait marquant : l'abondance d'articles publiés par Hermite sous forme de lettres, parmi lesquelles les lettres aux éditeurs des revues sont les plus nombreuses. Les noms de ces destinataires apparaissent donc dans notre base. Nous en avons déjà rencontré certains en discutant les noms de journaux, Borchardt, Liouville, par exemple. Mais d'autres s'y ajoutent, qu'Hermite n'utilise pas systématiquement pour désigner le journal dont ils s'occupent, mais à qui il adresse des lettres mathématiques pour publication : Francesco Brioschi pour les *Annali* ((2), 1859, p. 59 ou 2$^e$ s., (2), 1868, p. 97), Arthur Cayley pour les *Proceedings of the London Mathematical Society* ((4), 1873, p. 343-345), Paul Mansion pour la *Nouvelle Correspondance mathématique* ((2), 1876, p. 54-55), Gösta Mittag-Leffler pour *Acta mathematica* ((1), 1882, p. 368-370), James Joseph Sylvester ou Thomas Craig pour l'*American Journal of Mathematics* ((6), 1884, p. 173-175, (17), 1895, p. 6), Francisco Gomes Teixera pour le *Jornal de Ciências matemáticas e astronómicas* ((6), 1885, p. 81-82), Eduard Weyr pour *Casopispro pèstovàni mathematiky a fysiky* ((23), 1893, p. 273-274, (24), 1894, p. 65), et plusieurs autres. Hermite participe d'ailleurs plusieurs fois aux premiers numéros de journaux, ou de nouvelles séries de journaux.

Il faut remarquer que des lettres d'Hermite adressées à d'autres mathématiciens que les éditeurs sont aussi publiées directement : par exemple, des lettres de lui à Brioschi, Mittag-Leffler, Paul Gordan, Leo Königsberger, Hugo Gylden, Lazarus Fuchs, Ferdinand Lindemann, apparaissent dans le *Journal für die reine und angewandte Mathematik*, à côté de ses lettres à Borchardt. Le caractère de la « lettre à l'éditeur », genre qui peut être purement formel, s'estompe donc ici dans la mesure où nous savons qu'Hermite a entretenu des correspondances de longue durée avec certains de ces éditeurs-mathématiciens (comme Borchardt ou Sylvester, par exemple), comme avec d'autres destinataires de ses lettres publiées qui n'étaient pas des éditeurs (comme Königsberger). À l'inverse, ses échanges épistolaires avec d'autres mathématiciens, que ce soit ou non dans leur rôle d'éditeurs, semblent localisés à un résultat particulier, une occasion spécifique.

---

28  Si le nom historique de « Journal de Crelle » reste lui aussi en usage, c'est pour l'ensemble des numéros —ce qui constitue une différence essentielle avec l'usage d'Hermite. À ce propos, je n'ai pas utilisé la convention typographique courante de mettre en italique les noms de journaux, lorsqu'il s'agissait des dénominations d'Hermite. Ceci est une commodité personnelle ici ; il arrive que l'italique soit utilisé dans les *Œuvres* d'Hermite, pour *Journal de Crelle*, par exemple, mais de manière trop aléatoire pour en tirer des informations utiles (est-ce un choix des éditeurs, de la période, …?). Les questions typographiques restent un domaine mal exploré en histoire des mathématiques, voir toutefois les suggestions très intéressantes et variées de Reviel NETZ, *The Shaping of Deduction in Greek Mathematics*, Cambridge : Cambridge University Press, 1999 ; Jim RITTER, « Les nombres et l'écriture », in Yves MICHAUD (dir.), *Qu'est-ce que l'univers ?*, Paris: Odile Jacob, 2001, p. 114-129 ; Norbert VERDIER, *Joseph Liouville (1809-1882), son journal (1836-1874) et la presse de son temps : une entreprise éditoriale au service des mathématiques*, Thèse Université Paris-Sud, 2009, ch. II.

29  Charles HERMITE, *op. cit*. n. 23, vol. 4, p. 151 et 157.

30  Le cas du *Journal de mathématiques pures et appliquées* à l'époque de Liouville est traité dans Norbert VERDIER, *op. cit.* n. 28, à l'époque de Jordan dans Frédéric BRECHENMACHER, « Le Journal de M. Liouville sous la direction de Camille Jordan (1885-1922) », *Bulletin de la Sabix* (45), 2010, p. 65-71, mis en ligne le 09 octobre 2011, http://sabix.revues.org/730. Je renvoie à leur bibliographie pour les études, maintenant nombreuses, sur les (autres) journaux mathématiques.

Nous avons donc affaire à un continuum de relations professionnelles, sociales et mathématiques, les mêmes interlocuteurs pouvant intervenir comme responsables de journaux, correspondants, auteurs de résultats, sans qu'aucune de ces positions spécifiques en tant que telle ne fixe définitivement la nature, ou la profondeur, des liens à Hermite.

Une autre catégorie de personnes explicitée dans les articles renforce cette analyse, celle des traducteurs. Hermite précise en effet à plusieurs reprises qui a traduit tel ou tel article dont il se sert. Par exemple, dans une des lettres à Jacobi sur la théorie des nombres publiée en 1850, il indique : « C'est M. Borchardt lui-meme qui a bien voulu me traduire l'article de M. Gauss [sur les formes quadratiques ternaires][31] » ou encore, renvoie ainsi à un travail de Karl Weierstrass dans une note aux *Comptes rendus* de 1855 : « Voyez Journal de Crelle, tome 47, ou dans le Journal de Liouville (traduction de M. Woepcke), le Mémoire dans lequel ce savant géomètre a donné un aperçu de ses grandes et belles découvertes[32]. » Outre Borchardt et Franz Woepcke, on rencontre Victor Puiseux (comme traducteur de Jacobi), Irénée-Jules Bienaymé (comme traducteur de Pafnuti Lvovitch Tchebichef), ou Henri Padé (comme traducteur de Hermann Amandus Schwarz).

Comme on le voit, il ne s'agit pas de professionnels de la traduction, mais bien de mathématiciens utilisés comme traducteurs. Les exemples rapportés plus haut montrent d'ailleurs deux phénomènes a priori distincts. D'une part, les journaux français ne publient qu'en français (contrairement au journal de Crelle qui accepte des articles dans plusieurs langues), et n'hésitent pas à republier en traduction des articles parus ailleurs, en allemand et en russe en particulier ; on voit donc se mettre en place un service de traducteurs[33]. Ce qui frappe ici, c'est qu'Hermite mentionne ces traducteurs (et pas seulement l'auteur), lorsqu'il fait une référence explicite à l'article en traduction. D'autre part, le cas de Borchardt témoigne d'une aide plus personnelle et renvoie au fait bien connu qu'Hermite, grand admirateur des mathématiciens allemands, n'ait jamais maîtrisé correctement leur langue, malgré plusieurs séjours en Allemagne[34]. Bien que cela n'apparaisse pas à la lecture des articles, les deux aspects tendent à fusionner lorsqu'Hermite, au cœur de la vie académique dans la seconde moitié du siècle, utilisera ses réseaux variés pour promouvoir des traductions vers le français. Dans le cadre de la présente étude, nous pouvons seulement souligner encore une fois le mélange constant des liens unissant Hermite aux différentes catégories de personnes que ses œuvres mettent en scène : Borchardt est à la fois éditeur, correspondant, traducteur. Il est aussi auteur de résultats développés par Hermite. C'est vers ce dernier aspect que nous allons nous tourner maintenant.

Parmi les 144 personnes citées par Hermite et dont le décès est postérieur à 1800, 21 ne le sont que pour une des fonctions discutées plus haut, c'est-à-dire en tant que nom d'objet mathématique (« indicatrice de Dupin »), traducteur (Woepcke), éditeur ou destinataire de lettres en vue d'une publication dans un journal (Barnaba Tortolini, Ferdinand Caspary, Gergonne, Weyr, Craig, etc...). Cette situation touche d'ailleurs tant des personnes de générations antérieures à Hermite (Gergonne est né en 1771, un demi-siècle avant Hermite, Tortolini en 1808) que des mathématiciens beaucoup plus jeunes que lui (Weyr est né en 1852, Craig en 1855, Caspary en 1853, Gomes Teixera en 1851).

---

31     Charles HERMITE, *op. cit. n*. 23, vol. 1, p. 136.
32     Charles HERMITE, *op. cit. n*. 23, vol. 1, p. 473.
33     Cet aspect a été analysé pour le journal de Liouville par Norbert VERDIER, *op. cit. n*. 28, ch. 12 et 16 en particulier. Celui-ci met en évidence le rôle important joué par certains traducteurs, comme Olry Terquem ou Jules Houël.
34     Voir Thomas ARCHIBALD, « Charles Hermite and German mathematics in France » in Karen HUNGER PARSHALL and Adrian RICE (eds.), *Mathematics Unbound: the evolution of an international mathematical research community, 1800–1945*, History of Mathematics 23, Providence, RI : AMS, 2002, pp. 123–137 ; Catherine GOLDSTEIN, « Du Rhin et des nombres : quelques réflexions sur l'usage de la notion de transfert culturel en histoire des mathématiques », in Philippe HERT et Maurice PAUL-CAVALLIER (éd.), *Sciences et frontières*, Fernelmont : EME, 2007, p. 342-376.

Restent donc 123 noms propres, correspondant à des personnes au travail ou aux ouvrages mathématiques desquels Hermite se réfère. Comme c'est souvent le cas au 19$^e$ siècle[35], Hermite donne rarement des références bibliographiques complètes : apparaissent en général le nom de l'auteur et le titre du journal ou du livre (sous différentes formes, abrégées, par exemple), parfois la date ou le numéro du volume ; mais fréquemment, nous ne disposons que du nom de l'auteur et d'une indication plus ou moins précise du résultat visé. J'ai déjà évoqué l'emploi de l'avant-nom « M. ». Les noms sont parfois accompagnés (ou remplacés en deuxième mention) par d'autres avant-noms : « ce savant géomètre », « savant confrère », « l'illustre géomètre », « un géomètre italien distingué M. Betti », « l'un des meilleurs auteurs » (il s'agit de Christoph Gudermann), mais aussi de manière plus personnelle, « mon ami M. Sylvester » ou « un homme du mérite le plus distingué et donc la mémoire est restée chère à ses nombreux amis » (Olinde Rodrigues). Il arrive aussi que des précisions soient ajoutées, métier, lieu d'activité, par exemple : « M. Moutard, professeur à Paris», « M. Sylvester célèbre géomètre anglais», « M. G. Forestier, ingénieur des Ponts et chaussées à Rochefort », « M. Paul Gordan, professeur à l'université d'Erlangen ». Les mélanges existent, comme dans « un jeune géomètre du talent le plus distingué, M. Mittag-Leffler, professeur à l'université d'Helsingfors ». Ces mentions ne sont pas systématiques, elles semblent dépendre de la notoriété plus ou moins grande de la personne mentionnée par rapport au lieu de publication : les jeunes auteurs sont très souvent introduits, ainsi que les étrangers dans les textes destinés à une audience plus large, des étudiants par exemple. Leur effet est d'expliciter de manière particulièrement active une foule de relations individualisées, tant personnelles que professionnelles, dans les textes mathématiques d'Hermite. La grande majorité de ces 123 auteurs de travaux cités y apparaissent donc comme des personnes à part entière, et il semble donc pertinent de les étudier en tant que telles de plus près.

Un premier aspect est générationnel. Hermite étant né en 1822, j'ai distingué des périodes de 20 ans pour examiner cette question, l'une, baptisée par commodité dans ce qui suit « génération d'Hermite », correspondant aux années 1810-1830. L'ajustement est évidemment arbitraire. Ceci dit, ce découpage met Liouville, né en 1809 ou Gustav Lejeune Dirichlet, né en 1805, dans la génération précédente à celle d'Hermite, en accord avec le rôle qu'ils ont joué au début de sa carrière ; en revanche, Joseph Bertrand (né en 1822), Joseph Alfred Serret (né en 1819), et, dans un registre directement mathématique, Gotthold Eisenstein, Leopold Kronecker et Enrico Betti (nés en 1823), Francesco Brioschi (né en 1824), Arthur Cayley (né en 1821), Ernst Eduard Kummer (né en 1810) ou Karl Wilhelm Borchardt (né en 1817) se retrouvent avec cette normalisation dans la génération d'Hermite, une configuration qui a bien été pertinente dans le travail d'Hermite des années 1845-1860[36].

Les relevés obtenus sont consignés dans le tableau suivant : $n$ désigne la date de naissance, les lignes indiquent quelle période du travail d'Hermite est concernée (rappelons que la période de publication d'Hermite s'étend de 1842 à 1900). Par exemple (2$^e$ ligne, 2$^e$ colonne), entre 1860 et 1880, Hermite cite 17 auteurs nés entre 1790 et 1810.

---

35   Sur ce point, il existe de grandes variations selon les auteurs. Leur étude est un préliminaire pour toute enquête fondée sur les citations (en particulier celle décrite dans la première partie du présent article !), voir Catherine GOLDSTEIN, « Sur la question des méthodes... », article cité note 8.
36   Cette question joue un rôle dans la mise en place d'un champ de recherches que N. Schappacher et moi-même avons baptisé «arithmetic algebraic analysis », dans Catherine GOLDSTEIN et Norbert SCHAPPACHER, « A Book in Search of a Discipline (1801-1860) », in *The Shaping of Arithmetic, op. cit.* note 14, p. 3-65.

| Auteurs cités | n < 1790 | 1790 ≤ n < 1810 | 1810 ≤ n < 1830 | 1830 ≤ n < 1850 | 1850 ≤ n ≤ 1870 |
|---|---|---|---|---|---|
| avant 1860 | 5 | 11 | 22 | | |
| entre 1860 et 1880 | 10 | 17 | 30 | 16 | 2 |
| après 1880 | 8 | 12 | 22 | 20 | 15 |
| toutes dates | 12 | 24 | 39 | 31 | 17 |

    Les mêmes auteurs peuvent bien sûr être cités sur plusieurs périodes. Notons toutefois qu'aucune des listes de noms sous-jacentes à ces chiffres ne se stabilise, autrement dit que des auteurs cités à une période peuvent cesser ensuite de l'être, et d'autres apparaître, même dans les générations d'auteurs les plus anciennes. Parmi les 8 auteurs nés avant 1790 cités après 1880, figure par exemple Giovanni Plana, qui n'était pas apparu avant ; ceci est aussi vrai pour les auteurs nés entre 1790 et 1810, où figurent parmi les 12 Franz Neumann, Olinde Rodrigues ou Gabriel Gascheau, non cités auparavant, ou des auteurs de la génération d'Hermite (où n'apparaissent qu'après 1880 Pierre Alphonse Laurent ou Wilhelm Scheibner par exemple).

    Un autre aspect remarquable est la référence aux travaux de jeunes mathématiciens tout au long de la vie d'Hermite ; cet aspect est d'ailleurs souligné par le qualificatif « jeune » accompagnant parfois leur nom comme dans « jeune géomètre » que nous avons déjà rencontré pour Mittag-Leffler, mais qu'Hermite utilise aussi pour Émile Picard (né en 1856) par exemple, à la fin de la période 1860-1880. Il s'agit bien ici de lectures et d'utilisation de leurs travaux, et pas d'une simple connaissance ou d'un contact avec eux du fait même de l'évolution de la communauté scientifique, comme c'est le cas pour les destinataires de ses lettres qui sont éditeurs de journaux, comme George Chrystal, ou Craig évoqué plus haut.

    À l'inverse, pourtant, les années 1860 marquent une réelle coupure (ce sont d'ailleurs les années où Hermite obtient pour la première fois des positions stables). D'une part, seuls 38 auteurs différents sont cités pendant cette période, au lieu du double (75 et 78 respectivement) pour chacune des deux périodes suivantes. D'autre part, sur ces 38 auteurs, une vingtaine constituent le noyau stable des références d'Hermite, qui réapparaissent de manière récurrente dans son œuvre entière, toutes périodes confondues. Il s'agit de Joseph-Louis Lagrange, Adrien-Marie Legendre, Carl-Friedrich Gauss, Augustin Louis Cauchy parmi ceux nés avant 1790, de Niels Abel, Dirichlet, Jacobi, Liouville, Charles Sturm, nés entre 1790 et 1810, et dans la génération d'Hermite même, de Borchardt, Brioschi, Charles Auguste Briot, Cayley, Eisenstein, Évariste Galois, Kronecker, Joseph-Alfred Serret, Sylvester, Tchebychef et Weierstrass. Leur rapprochement témoigne de prédilections thématiques récurrentes autour de la théorie des nombres, des équations, des séries et des fonctions elliptiques. Le recours à certains d'entre eux est intensif, avec des références fréquentes et multiples : Jacobi est ainsi cité dans plus de 50 articles d'Hermite, Gauss dans 27, Lagrange, Cauchy et Kronecker dans 24.

    Les autres références sont plus ponctuelles, et après 1860, la diversification des noms va de pair avec un usage localisé, pour quelques articles, alors que le nombre de noms cités par article est en croissance nette. Parmi les 18 nouveaux auteurs nés après 1830 qui apparaissent dans les références postérieures à 1860, seuls 5 (Lazarus Fuchs, Hugo Gylden, Rudolf Lipschitz, Gösta Mittag-Leffler et Edmond Laguerre) sont encore mentionnés dans les articles d'Hermite après 1880 ; le renouvellement des sources est donc important.

    Comme les noms déjà donnés l'ont suggéré, un grand nombre de nationalités et de régions d'exercice sont représentées : France, états allemands et italiens, Bohème, Suisse, Danemark,

Norvège et Suède, Finlande, Russie, Pays-Bas, Belgique, Royaume-Uni[37]. Mais les proportions sont très différentes : 47 auteurs référencés travaillent en France (auxquels s'ajoute Thomas Stieltjes, originaire des Pays-Bas, à partir du milieu des années 1880), 40 en Allemagne, 9 en Italie et en Russie, 5 au Royaume-Uni (dont Sylvester qui passe aux États-Unis une partie de la période), 3 en Suède, 2 en Belgique, en Finlande et au Danemark, 1 en Bohème et en Suisse. En l'absence d'études quantitatives sur la scène européenne de cette époque, il est bien sûr difficile d'évaluer la singularité de cette répartition. Mais la capacité d'Hermite à servir d'intermédiaire aux mathématiciens allemands est ici bien mise en évidence[38].

Il n'y a pas de différence frappante selon les périodes, on retrouve pour toutes une nette domination franco-allemande des références. La liste des auteurs cités dans toutes les périodes (donnée plus haut) est d'ailleurs en parfaite cohérence avec la répartition globale : 7 Allemands, 8 Français, 2 Anglais, 1 Norvégien, 1 Russe, 1 Italien. En revanche, croiser cette information avec celle des générations est éclairant. La plupart des auteurs cités nés avant 1790 sont Français (l'exception notable étant Gauss). La répartition entre Français et Allemands est équilibrée pour la génération précédant celle d'Hermite, puis bascule en faveur des mathématiciens allemands pour les auteurs référencés nés entre 1810 et 1850 : par exemple, Hermite cite 10 auteurs allemands de sa génération contre 6 Français avant 1860, 12 contre 9 entre 1860 et 1880, 9 contre 7 après 1880. En revanche, 6 mathématiciens travaillant en France et nés après 1850 sont cités à partir des années 1880 (pour 2 seulement travaillant en Allemagne, 2 en Italie, 2 en Finlande, 1 en Bohème) : on assiste ici à l'entrée de Paul Appell, de Henri Poincaré, de Henri Padé, etc.

La majorité des auteurs cités sont professeurs dans une université, un institut technique, une école d'ingénieurs (Bauschule, école centrale) ou travaillent dans un observatoire. Les exceptions sont presque exclusivement françaises : professeurs de mathématiques en lycée, ils sont le plus souvent issus de Polytechnique, y étant parfois répétiteurs ou chargés de cours (Théodore Florentin Moutard), ou sont passés par l'École normale supérieure (Charles Joubert, Edouard Lucas). Les autres (Louis Bourguet, Charles Biehler, …) ont fait des thèses avec Hermite. Plus d'une trentaine des mathématiciens français cités sont élus à l'Académie des sciences. Il s'agit donc d'un monde bien délimité, les travaux cités ayant été publiés dans des revues de mathématiques professionnelles usuelles.

---

37    Je ne tiendrai pas compte dans cette brève présentation des variations politiques pendant la vie d'Hermite et parlerai par exemple librement d'Allemagne et d'Italie, d'autant plus que la circulation des mathématiciens entre leurs postes successifs dans ces deux régions rend assez vaine une classification plus fine.
38    Voir la note 34 ci-dessus. Notons tout de même qu'Hermite publie beaucoup dans des journaux étrangers, et qu'une étude des références utilisées selon le journal de publication, étude que je ne présenterai pas ici, serait nécessaire sur ce point.

| Journal | Date | Nom | Contexte |
|---|---|---|---|
| Proceedings of the London Mathematical Society, 7, p. 173-175 | 1876 (lu le 13 avril 1876) | H.J.S. Smith | ajoute une note bibliographique avec la permission d'H. |
| | | M. Heine | H. simplifie ici sa méthode pour prouver un théorème d'Eisenstein |
| | | Eisenstein | |
| | | Crelle | journal de |
| | | Kepler | équation de |
| | | M. Serret | |
| Journal für die reine und angewandte Mathematik, 81, p. 220-228 | 1876 | M. L. Königsberger | l'article est extrait d'une lettre de H. à lui du 2 octobre 1875 |
| | | Maclaurin | formule de |
| | | Bernoulli | nombre de |
| | | Cauchy | théorème de |
| | | M. C.-O. Meyer | |
| | | M. Hugo Gylden | |
| Nouvelle Correspondance mathématique, 2, p. 54-55 | 1876 | M. Paul Mansion | l'article est extrait d'une lettre de H. à lui du 25 novembre 1875 |
| | | M. Delaunay | |
| | | Liouville | journal de |
| Annales de la société scientifique de Bruxelles, I, p. 1-16 | 1876 | Crelle | journal de |
| | | Legendre | |
| | | Jacobi | Sa généralisation d'un résultat de Legendre est le point de départ d'H |
| | | Clebsch | |
| | | Göpel | |
| | | M. Rosenhain | |

Illustration 2: Noms propres dans les articles de C. Hermite (extrait de la base de données)

III Retour sur Charles Hermite

La mise en relation de ces deux enquêtes conduit à un constat intriguant : Hermite cite peu, voire pas du tout, les auteurs du réseau de textes H-K, c'est-à-dire celui pour lequel il est un des pères fondateurs. Ceci est d'autant plus frappant qu'il connaît leurs travaux (ne serait-ce que parce qu'il les présente à l'Académie pour publications dans les *Comptes rendus*). Il les encourage même dan certaines directions, comme nous l'avons vu pour Poincaré un peu plus haut, et comme ses correspondances en témoignent abondamment. Il reconnaît d'ailleurs leur filiation, par exemple dans son intervention lors de son propre Jubilé en 1892 : « M. Camille Jordan a été bien au-delà de mes premières tentatives, dans la question arithmétique de la réduction des formes, dans l'étude des équations algébriques et la théorie des substitutions qui en est le fondement[39] », écrit-il alors que Jordan n'apparaît dans notre deuxième enquête que comme éditeur du *Journal de mathématiques*. Ce constat s'étend d'ailleurs aux auteurs étrangers, comme Paul Gordan, qui poursuit le travail sur les invariants ou Hermann Minkowski, dont la *Geometrie der Zahlen* est dédiée à Hermite et qui est considéré par plusieurs mathématiciens français comme le « continuateur d'Hermite »[40]. Les apports des deux prosopographies ne sont donc pas redondants, ils éclairent deux aspects distincts de l'insertion d'Hermite dans le monde mathématique de la deuxième moitié du $19^e$ siècle.

Les deux enquêtes sont d'ailleurs loin d'épuiser la richesse de ce monde. Hermite, nous l'avons vu plus haut, a entretenu une vaste correspondance, dont une infime partie a été entrevue ici sous forme de lettres publiées dans les journaux réguliers. Or, ces correspondances font apparaître des ensembles de noms souvent très différents de ceux mentionnés dans les articles. Pour ne prendre qu'un exemple, les lettres d'Hermite à Thomas Stieltjes mettent en scène 220 personnes environ (contre rappelons-le, 159 seulement dans les articles de recherche)[41]. Certains noms de collègues comme Benjamin Baillaud, Désiré André, David Bierens de Haan ou Sofia Kowalevskaia (« Mme de Kowalewski ») font leur apparition alors qu'ils ne sont pas sollicités dans le cadre des articles. Il s'agit bien d'une ouverture sur de nouvelles relations ou de nouvelles activités, qui participeraient pleinement de l'élaboration d'une biographie scientifique, et ne laissent presque pas de traces explicites dans les articles de recherche. Hermite écrit ainsi en 1883 que : « M. Tisserand que j'ai eu pour élève est un de mes amis et j'ai avec tous [les astronomes de l'Observatoire de Paris] de bons et excellents rapports » ou encore en 1889 que c'est «…Mr Darboux avec qui je m'entretiens de vos intérêts. C'est ordinairement le lundi que j'ai l'occasion de le voir et de causer amicalement avec lui dans le cabinet de M. Bertrand »[42].

De nouvelles catégories de personnes sont aussi explicitées par Hermite, par exemple sa famille. On rencontre « [s]a famille de Lorraine », « [m]a mère », « Mme Hermite », etc., qui ne jouent pas de rôle dans ses articles, mais aussi en qualité familiale, plusieurs noms déjà rencontrés dans les deux prosopographies : le beau-frère d'Hermite, Joseph Bertrand, ses deux gendres Georges Forestier et Émile Picard, son neveu par alliance Paul Appell. Les deux premiers sont polytechniciens, les deux derniers normaliens, et si rien ne les distingue dans les articles des autres auteurs cités, on constate dans les correspondances que les réunions familiales sont l'occasion de discussions scientifiques. La biographie collective (familiale) que cette situation suggère rejoindrait d'ailleurs une des sources du genre prosopographique, la généalogie familiale des élites[43].

---

39  Charles HERMITE, *op. cit. n*. 23, vol. 4, p. 583. Sur les liens entre Jordan et Hermite, voir Frédéric BRECHENMACHER, « Self-portraits with Évariste Galois (and the shadow of Camille Jordan)», *Revue d'histoire des mathématiques* (17), 2011, p. 271-369.
40  Sur les recherches de Gordan sur la théorie des invariants, voir Karen HUNGER PARSHALL, « Toward a History of Nineteenth-Century Invariant Theory », in David ROWE and John McCLEARY (eds.), *The History of Modern Mathematics, vol. I = Ideas and their Reception*, Boston, … : Academic Press, 1989, p. 157-206 ; sur les relations de Minkowski à Hermite et ses successeurs, voir Sébastien GAUTHIER, *La Géométrie des nombres comme discipline*, Thèse de l'UPMC, Paris, 2007 et « Justifier l'utilisation de la géométrie en théorie des nombres : des exemples chez C.F. Gauss et H. Minkowski », in Dominique FLAMENT et Philippe NABONNAND, *Justifier en mathématiques*, Paris : MSH, 2011, p. 103-128.
41  Charles HERMITE et Thomas STIELTJES, *Correspondance*, éditée par Benjamin BAILLAUD et Henry BOURGET, 2 vols. Paris: Gauthier-Villars, 1905. Cette édition contient un index d'auteurs.
42  Charles HERMITE et Thomas STIELTJES, *op. cit.* n. 41, vol. 1, p. 12 et p. 374 resp.

La catégorie familiale est donc tout aussi fluide que les autres catégories explicitées par Hermite (traducteurs, éditeurs de journal) et rencontrées plus haut. Un exemple supplémentaire témoigne de l'imbrication constante des rôles, et est particulièrement révélateur des difficultés d'interprétation que toute entreprise de catégorisation suscite, même lorsqu'elle s'appuie sur les termes mêmes des acteurs : c'est la catégorie d'élèves. Hermite utilise ce qualificatif à de nombreuses reprises. En réunissant articles de recherche et correspondance avec Stieltjes, par exemple, nous voyons ainsi apparaître : « un de mes élèves M. Bourguet », « un de mes élèves, M. Poincaré, ingénieur des Mines, professeur à la Faculté des sciences de Caen », « le doyen de la Faculté de Toulouse, M. Baillaud, qui est un de mes élèves, est aussi directeur de l'Observatoire », « M. Désiré André, un de mes élèves et mathématicien distingué, …», « M. Tisserand que j'ai eu pour élève...», « un de mes élèves qui est est aussi l'élève de M. Weierstrass, M. Mittag-Leffler », etc. Certaines des personnes désignées ont fait une thèse suivie de près ou de loin par Hermite, mais ce n'est pas le cas pour tous. Avoir suivi un cours d'Hermite n'est pas non plus une bonne définition pour la catégorie : Jules Maillard de la Gournerie qui a rédigé des cours d'Hermite n'est pas mentionné comme tel. Comme d'autres, le mot semble surtout tisser, voire raviver des liens plus personnels, et sur un plus long terme : une fonction a priori professionnelle est ici comme détournée pour suggérer du privé au sein de l'institutionnel et des positions officielles des protagonistes.

L'histoire des mathématiques, ou des sciences en général, est souvent le moyen de témoigner des aspects humains de la science. La biographie participe pleinement de cette problématique : montrer du héros les aspects singuliers et touchants, au sens le plus large. La biographie peut aussi tenter de résoudre les conflits et les oppositions qui ont traversé le domaine dans les dernières décennies—comme cela été le cas en histoire générale : ceux entre stabilité des résultats et actions pour les obtenir, entre les grandes masses et les micro-mouvements, l'avènement d'un nouveau point de vue sur le monde, d'une nouvelle manière de penser ou de prouver, et les hésitations ou les mesquineries du quotidien de la science. Dans la biographie, ces différentes échelles peuvent s'articuler, ces différentes oppositions peuvent se trouver apparemment résolues. Le savant, perçu comme être individuel, est restitué dans un milieu, un état de la science de son époque, bref un être social qui participe de l'élaboration de la science, bien collectif par exemple. La biographie comme genre ne s'oppose alors pas tant à l'histoire conceptuelle ou à celle des institutions qu'elle ne sert à les réconcilier. Au début de sa biographie de James Joseph Sylvester, Karen Parshall écrit de manière caractéristique : «This book aims to tell, for the first time, the complex story of Sylvester's life by situating that life as fully as possible within the political, religious, mathematical and social currents of nineteenth-century England. It aims to demythologize the man by placing him in his milieux at the same time it demystifies his mathematics by revealing it as a very human endeavour.…The story of Sylvester's life, for example, also sheds light on the evolution of mathematical thought, the ways in which mathematics may be done, and what factors may shape the mathematician's ideas.…and it highlights the very human side —shaped by factors as diverse as religion, ego and depression — of what many view as that most inhuman and otherwordly of intellectual endeavors : mathematics[44]. »

Nous avons proposé ici une voie quelque peu différente, fondée sur d'autres convictions sur l'activité mathématique. Hermite y apparaît néanmoins dans ses singularités : sa fidélité à une poignée d'auteurs et à leurs ouvrages, Jacobi, Gauss et parmi ses contemporains Kronecker, se couple avec une ouverture aux travaux de jeunes auteurs et une variation progressive de ses thématiques de recherche qui, du même coup, ne l'attache pas, mathématiquement du moins, à la discipline de recherches que ses propres travaux ont contribué à créer en France. Sa situation internationale a pu aussi être précisée : une attention privilégiée aux mathématiciens allemands, en particulier à travers une lecture régulière du *Journal für die reine und angewandte Mathematik*, ne l'empêche pas d'accompagner l'ouverture de la scène internationale. D'une part, il répond volontiers, surtout dans ses dernières décennies, aux sollicitations de collègues qui lancent ou reprennent des revues mathématiques dans des pays variés, en leur envoyant des articles ou des lettres mathématiques, du *Jornal de sciencias mathematicas e astronomicas*, de Gomes Texeira au Portugal au *Casopispro pèstovàni mathematiky a fysiky* de Weyr à Prague ou l'*American Journal of Mathematics*, lorsqu'il est repris par Craig. D'autre part, réciproquement, il contribue à la diffusion en France de travaux variés, soit en encourageant des traductions, soit en servant d'intermédiaire pour les *Comptes rendus* de l'Académie des sciences. À une période de tensions religieuses et de rivalités entre institutions, ceux qu'il cite avec des mentions favorables (« mon ami », « homme de mérite », etc.) incluent des polytechniciens et des normaliens, des catholiques[45], des juifs, des incroyants déclarés, des royalistes et des républicains —et ce malgré les convictions personnelles tranchées (catholique anti-Républicain) que ses lettres font apparaître par ailleurs. Les rôles qu'assignent ses descriptions sont provisoires et entremêlés, tissant sans cesse des continuités entre de multiples activités du monde mathématique de son temps. Si j'ai eu recours à des prosopographies pour cerner ces aspects, Hermite n'est donc pas ici un membre d'un collectif, ordinaire ou non, ni un individu situé *au milieu de* contextes scientifiques ou sociaux. Il est bien, du point de vue de l'historiographie, examiné comme un (ou des) collectifs mêmes : que nous soient accessibles à l'intérieur de ce collectif certains déplacements, certains changements, par exemple dans les références mobilisées par Hermite ou dans les catégories descriptives dont il joue, témoigne qu'une telle représentation n'exclut pas de restituer également la liberté d'action de la personne en tant que telle.

---

43    Cette généalogie et l'intrication des relations familiales et professionnelles de Joseph Bertrand, donc d'Hermite, sont explorées dans Martin ZERNER, « Le règne de Joseph Bertrand (1874-1900) », in Hélène GISPERT, *op. cit.* n. 12 , p. 298-322. Sur les liens avec l'histoire de la prosopographie, voir Karl Ferdinand WERNER, « L'apport de la prosopographie à l'histoire sociale des élites », in Katharine S. B. KEATS-ROHAN, ed., *Family Trees and the Roots of Politics*, Woodbridge : Boydell, 1997, p. 1-21.
44    Karen HUNGER PARSHALL, *James Joseph Sylvester: Jewish mathematician in a Victorian world*, Baltimore: Johns Hopkins University Press, 2006.

---

45 Plusieurs religieux (Charles Joubert, Charles Biehler, ...) apparaissent dans les auteurs cités et une enquête sur les réseaux catholiques d'Hermite donnera sans doute des résultats particuliers, Hermite étant impliqué par exemple dans l'ouverture de l'Institut catholique de Paris. Ceci n'aboutit à aucune exclusion dans les références d'articles (ni dans les soutiens qu'il accorde, comme en témoigne sa correspondance).